\documentclass[12pt]{amsart}
\usepackage[latin1]{inputenc}
\usepackage{amsmath}
\usepackage{amsfonts}
\usepackage{amssymb}

\usepackage{fullpage}

\newtheorem{lem}{Lemma}
\newtheorem{prop}{Proposition}
\newtheorem{thm}{Theorem}
\newtheorem{cor}{Corollary}
\theoremstyle{definition}
\newtheorem{defn}{Definition}
\newtheorem{ex}{Example}

\newcommand{\Z}{\mathbb{Z}}
\newcommand{\R}{\mathbb{R}}

\newcommand{\CP}{\mathbb{CP}}

\newcommand{\bp}{\bar{p}}
\newcommand{\bq}{\bar{q}}
\newcommand{\bE}{\bar{E}}

\newcommand{\bv}{\bar{v}}
\newcommand{\bg}{\bar{\gamma}}
\newcommand{\bal}{\bar{\alpha}}

\newcommand{\al}{\alpha}
\newcommand{\g}{\gamma}

\newcommand{\inv}{^{-1}}

\newcommand{\pr}{\partial_r}
\newcommand{\pth}{\partial_{\theta}}

\newcommand{\sig}{\Sigma}

\newcommand{\sigE}{\Sigma_pE}

\renewcommand{\prod}[2]{\langle #1 ,#2\rangle}

\newcommand{\ip}[2]{\ensuremath{\langle #1,#2\rangle}}
\newcommand{\bfu}{\mathbf{u}}
\newcommand{\bfv}{\mathbf{v}}

\DeclareMathOperator{\dis}{d}

\begin{document}

\title{Submetries vs. submersions}

\author{Luis Guijarro}
\address{ Department of Mathematics, Universidad Aut\'onoma de Madrid, and ICMAT CSIC-UAM-UCM-UC3M}
\curraddr{}
\email{luis.guijarro@uam.es}
\thanks{The first author was supported by research grant MTM2008-02676-MCINN}

\author{Gerard Walschap}
\address{Department of Mathematics, University of Oklahoma.}
\curraddr{}
\email{gerard@ou.edu}
\thanks{}

\maketitle

\section{Introduction and main results}

Riemannian submersions are one of the main tools used in Riemannian geometry. On the one hand, they are a basic component of the structure of certain metrics: for example, the metric projection of an open manifold with nonnegative sectional curvature onto its soul is a Riemannian submersion. On the other, every known example of positive curvature is constructed by means of these maps.

In \cite{bere}, Berestovskii introduced a purely metric version of  Riemannian submersions: a map $\pi:X\to B$ between metric spaces is a \emph{submetry} if for every $p\in X$, any closed ball $B(p,r)$ of radius $r$ centered at $p$ maps onto the ball $B(\pi(p),r)$.    
The increasing use of Gromov-Hausdorff distance in Riemannian geometry is lending growing relevance to this concept. In \cite{BG} it was proved that submetries between Riemannian manifolds are $C^{1,1}$ Riemannian submersions; thus, new features of submetries appear only under the presence of some type of singularities, whether in the base or the domain of $\pi$.

This note has several purposes: first, we  give examples that show how some of the well known splitting theorems for Riemannian submersions fail in the context of submetries (recall that a submersion $\pi:M\to B$ is said to \emph{split} if $M$ is locally a metric product and $\pi$ is projection onto one of the factors). Among these, we exhibit a submetry from a sphere with a locally flat metric everywhere except for a codimension one singular set (in contrast, a Riemannian submersion from a compact flat manifold always splits). This  motivates a further study, in sections \ref{quasigeodesics} and \ref{extremals}, of how metric singularities interfere with the fiber structure of submetries. One highlight of the latter section is the fact that extremal sets cannot be oblique with respect to the horizontal-vertical structure that a submetry introduces in the total space. Finally, we turn to 
some splitting theorems. 

Throughout the paper, Alexandrov spaces will be finite dimensional, and $|pq|$ stands for the distance $\dis(p,q)$ between the points $p$ and $q$. 
For convenience of the reader, we have included a preliminary section with some of the main results in \cite{Lyt} that will be used in the rest of the paper.

This work was initiated during the workshop "Manifolds with nonnegative sectional curvature" that was held at the American Institute of Mathematics on September 2007. The authors want to thank the Institute for the excellent working conditions that made this research possible. We would also like to thank A. Lytchak for his constructive criticism. 

\section{Structure of submetries between Alexandrov spaces}

In this section we collect several results that will be needed in what follows, starting with lifts. Recall that if $\bg:[0,a]\to B$ is a geodesic (i.e, shortest curve), a geodesic lift $\g:[0,a]\to X$ is a geodesic such that $\pi\circ\g=\bg$. Its existence and uniqueness has been proved in for instance, \cite{BG}, although only for the case in which $\bg$ is the unique geodesic between its endpoints. It is, however, easily shown that this latter condition is not necessary, although at the cost of losing uniqueness of the lift; a more general version of this lemma appears as
Lemma 4.4 in \cite{Lyt_open}, but we include one here for the reader's convenience. 

\begin{lem}[Existence of geodesic lifts]\label{lifts}
Let $\pi:X \to B$ denote a submetry between finite dimensional Alexandrov spaces. For any unit speed geodesic $\bg:[0,a]\to B$, and any $p\in\pi\inv(\bg(0))$, there is a geodesic lift $\g:[0,a]\to X$ of  $\bg$ starting at $p$, with length equal to that of $\bg$. 
\end{lem}
\begin{proof}
We will construct $\g$ as a limit of piecewise geodesics. For any integer $n\geq 0$, choose a subdivision of $[0,a]$ with points $t_i^n=i\cdot a/n$ and denote by $\bp_i^n$ the points $\bg(t_i^n)$. 
Construct a sequence of points $\{p_i^n\}$ in $X$ by letting $p_0^n=p$ and $p_i^n$ be a point in the fiber over $\bp_i^n$  such that $|p_{i-1}^n p_i^n|$ coincides with $|\bp_{i-1}^n \bp_i^n|$; let $\g^n$ be the curve obtained by connecting each $p_i^n$ with its sucessor $p_{i+1}^n$. Ascoli-Arzela's theorem ensures the existence of some subsequence converging to a curve $\g:[0,a]\to X$, which is clearly a lift of $\bg$. Because of the semicontinuity of length under limits of curves, and the fact that $\pi$ is a submetry, we obtain that the length of $\g$ agrees with that of $\bg$.
\end{proof}
It is also clear that the lift constructed is \emph{horizontal}, meaning that its direction at each time lies in the sets $H_{\g(t)}$ defined later in this section.

Non-uniqueness of liftings occurs for instance in the canonical submetry that sends the double disk $D(\Delta)$ to the disk $\Delta$: any geodesic connecting two points in the boundary of $\Delta$ can be lifted in two different ways to $D(\Delta)$. 

Next, we recall some results of Lytchak about the structure of general submetries between Alexandrov spaces \cite{Lyt}. 
Denote by $\pi:X\to B$ a submetry between two Alexandrov spaces. We will assume all along that $X$, $B$ are complete, although some of the statements in this section apply locally. Unless stated otherwise, the dimensions of $X$ and $B$ are $n+k$ and $n$ respectively.

%
%
Alexandrov geometry provides differentiable tools for further study of submetries. The term ``differentiable'' must, of course, be understood in a metrical sense:

\begin{defn} Let $f:X\to Y$ a map between Alexandrov spaces.
$f$ is said to be \emph{differentiable at $p\in X$} if  the maps $f^t:(tX,p)\to (tY,f(p))$ converge, in the Gromov-Hausdorff sense,  to some limit map $D_pf:T_pX\to T_{f(p)}Y$ as $t\to\infty$.
 In this case, the limit map $D_pf$
 is called \emph{the derivative of $f$ at $p$}.\end{defn}

Reference \cite{Lyt} deals with more general spaces, but in the Alexandrov case the above definition suffices. Observe also that in this situation, $T_pX$ is actually the Euclidean cone $C(\Sigma_p)$ over the space $\Sigma_p$ of directions at $p$, and similarly for $T_{f(p)}Y$.

\begin{defn}
A \emph{homogeneous submetry} between Euclidean cones is a submetry $g:C(\sig)\to C(S)$ commuting with dilations; i.e, $g(tp)=tg(p)$, $t\in\R$.
\end{defn}

The following key feature of the differentials of submetries is found in \cite{Lyt}: 
\begin{prop}
Submetries are differentiable at any point $p\in X$, and its derivative $D_p\pi:T_pX\to T_{\pi(p)}B$ is a homogeneous submetry.
\end{prop}

In light of this property, the following result of Lytchak's thesis \cite{Lyt} is of paramount importance, as it mimics the usual splitting of tangent spaces into horizontal and vertical subspaces:

\begin{prop}\label{hover}
Let $\sig$, $S$ denote two Alexandrov spaces with curvature $\geq 1$. If $f:C(\sig)\to C(S)$ is a  homogeneous submetry between their Euclidean cones, then 
\begin{enumerate}
\item The fiber over $0$ is a totally convex subcone $C(V)$ of $C(\sig)$ where $V$ is a totally convex subset of $\sig$. 
\item The cone $C(H)$ over the polar set $H\subset\sig$ to $V$ (that is, the set of points in $\sig$ at distance at least $\pi/2$ from V) agrees with the cone over the set of horizontal directions; i.e, those directions $h\in \sig$ with $|f(h)|=|h|$. Furthermore, $H$ is a totally convex subset of $\Sigma$.
\item $\dis(V,H)=\pi/2$, and for any $x\in\Sigma\setminus(V\cup H)$ there is precisely one minimal geodesic between $H$ and $V$ that passes through $x$.
\end{enumerate} 
\end{prop}

We remind the reader that in Alexandrov geometry, a totally convex set is one which contains all minimal geodesics whose endpoints lie in the set. An easy mistake to make would be to assume that the space of directions at some $p\in X$ is necessarily the spherical join of the horizontal and the vertical spaces. This is not the case, as the following example shows:

\begin{ex}
Let $X$ be the Euclidean cone over $\CP^2$ with its Fubini-Study metric, and choose any ray from the vertex, $\g:[0,\infty)\to X$. Its associated Busemann function $b_{\gamma}$ is a submetry onto $[0,\infty)$.  At the vertex, the horizontal space is a $\CP^1$, while its vertical space is a point. However its space of directions is the whole $\CP^2$ (see also section (6.4) in \cite{Pe}). 
\end{ex}

In order to overcome this difficulty, the following result from Lytchak is needed:

\begin{prop}
Let $f:C(\sig)\to C(S)$ be a homogeneous submetry, $x\in C(\sig)\setminus \big(C(V)\cup C(H)\big)$. Then there is a unique pair $h\in C(H)$, $v\in C(V)$ with $x=h+v$ and $\langle h,v\rangle=0$. In fact such pair is formed by the projections of $x$ on $C(H)$ and $C(V)$; furthermore $f(x)=f(h)$. 
\end{prop}

The next construction, which appears at several points in \cite{Lyt}, will be needed to describe extremal sets:

\begin{defn}
Let $A\subset H$, $B\subset V$ two sets. We denote by $P(A,B)$ the set of all points that belong to the image of geodesics connecting points of $A$ to points of $B$, with the convention $P(A,\emptyset)=A$, $P(\emptyset,B)=B$.
\end{defn}

The following facts are straightforward:
\begin{itemize}
\item If $A$ and $B$ are totally geodesic in $H$ and $V$, then $P(A,B)$ is totally geodesic in the space of directions $\sig_p$;
\item $P(H,V)=\sig_p$.
\end{itemize} 

We finish this section wtih a few minor observations:

\begin{cor}
Let $\pi:X \to B$ a submetry with connected fibers. If $p\in X$ has space of directions of diameter less than $\pi/2$, then the fiber through $p$ coincides with $p$.
\end{cor}

\begin{proof}
This is an immediate consequence of Proposition \ref{hover}.
\end{proof}

It should be observed that the situation of the corollary can occur in practice; for instance, let $B$ be the spherical suspension of $\CP^n$, and $f:S^{2n+1}\to\CP^n$ the Hopf map. Then the suspension of $f$ gives a submetry $\bar{f}:S^{2n+2}\to B$ whose fiber over a vertex point of $B$ is a ``vertex'' point of $S^{2n+2}$ when seen as a suspension of $S^{2n+1}$.   

\begin{cor}
If $b\in B$ is a regular point of $B$ (i.e, if the space of directions $\sig_b$ is isometric to the unit sphere $S^{n-1}$), then for any $p\in F_b$, we have a metric splitting 
$\sig_p= S^{n-1}\ast A_p$, for a positively curved Alexandrov space $A_p$ formed by the unit vertical directions at $p$. 
\end{cor}

\begin{proof}
This follows from Proposition \ref{hover} and the well known fact that any Alexandrov space with both curvature and diameter $\geq 1$ is a spherical suspension (see \cite{BBI} for instance).
\end{proof}


\section{examples}

In this section, we show by means of examples that some of the well known facts concerning Riemannian submersions do not hold for submetries. 
\subsection{The ``flat'' Hopf fibration}

The usual Hopf fibration $S^3\to S^2$ is a  Riemannian submersion when $S^3$ and $S^2$ are given the round metrics of constant curvature 1 and 4 respectively. We will modify the metric on $S^3$ to obtain an Alexandrov metric with the following characteristics:
\begin{itemize}
\item The metric is locally isometric to a flat metric everywhere except in a hypersurface, where it lacks smoothness;
\item It projects to a metric in $S^2$ for which the Hopf map $f:S^3\to S^2$ is a submetry.
\end{itemize} 

Recall that the sphere $S^3$ can be written as the union of two solid tori: the set of points with $x_0^2+x_1^2\geq 1/2$, and those with $x_0^2+x_1^2\leq 1/2$. Clearly each one of these tori is a union of Hopf fibers.  Identifying each tori with 
$S^1\times D^2$, where $D^2$ is the unit disc in the plane $\R^2$, the boundary identification is given by $(x_0,x_1,x_2,x_3)\mapsto (x_2,x_3,x_0,x_1)$ along the Clifford torus. 

To construct the flat metric, endow each torus $S^1\times D^2$  with the standard flat product metric. This yields two nonnegatively curved Alexandrov spaces with isometric boundaries under the above identification. Thus, by Petrunin's gluing theorem \cite{Pet}, the metric on $S^3$ obtained by gluing the boundaries is an Alexandrov metric with nonnegative curvature. It is clear that the Hopf action is still by isometries, thus providing a submetry $\pi:S^3\to S^2$ as claimed. 

It is interesting to observe that nonetheless, along the smooth part of the submetry, the usual O'Neill's $A$-tensor does not vanish: if $X$ is the unit vector tangent to the $S^1$ factor, and $\pr$, $\pth$ are the polar coordinate vectors in $D^2$, then the Hopf fiber is tangent to $X+r\pth$, and its orthogonal complement is spanned by $Y_1=rX-\pth$, $Y_2=\pr$. An easy computation shows that 
$$
A_{Y_1}Y_2=\frac{1}{2} \frac{\prod{-X}{X+r\pth}}{1+r^2}\left(X+r\pth\right)= -\frac{1}{2(1+r^2)}\left(X+r\pth\right).
$$

In the absence of singularities, a flat metric in a compact total space implies flatness of the base and global vanishing of the $A$ tensor (see for example \cite{GP} and \cite{Wa}). However, the presence of the metric singularity along the Clifford torus excludes this behavior. 


\subsection{Totally geodesic fibers are not always isometric} 
Given a Riemannian submersion, recall that any curve in the base space induces a diffeomorphism between the fibers over the endpoints, obtained by assigning to each point $p$ in the initial fiber the endpoint of the horizontal lift through $p$ of the curve. These maps are called \emph{holonomy diffeomorphisms}. 
If the fibers of a Riemannian submersion are totally geodesic, then the holonomy diffeomorphisms between fibers are isometries. This is no longer true for submetries, as our next example shows: 

Consider the action of $\Z$ on $\R^2\times \R=\mathbb{C}\times\R$ given by 
\[m(z,t) =(e^\frac{im\pi}2z, t+m),\qquad m\in\Z,\quad z\in\mathbb{C}, \quad t\in\R.
\]
 Let $M$ denote the quotient space by this action; it clearly inherits a flat metric. The action of $\R$ by translations on the second factor of $\mathbb{C}\times\R$ commutes with that of $\Z$, and hence induces an isometric action on $M$. The quotient space $B$ is isometric to a euclidean cone over a circle of length $\pi/2$, and there is a submetry $\pi:M\to B$. The fibers of $\pi$ are totally geodesic in $M$; however, they are not isometric: the one corresponding to the vertex of the cone has length 1, while all the others have length $4$. These statements are easier to visualize by observing that $M$ is just the mapping cylinder of the rotation of the plane by an angle $\pi/2$ with its flat metric, while the fibers of $\pi$ correspond to the subsets obtained from the vertical factor.

\section{Quasigeodesics and submetries}\label{quasigeodesics}

Quasigeodesics are natural substitutes for geodesics in Alexandrov spaces. Our aim in this section is to determine the extent to which the behaviour of horizontal geodesics in submersions carries over to quasigeodesics in submetries. The reader can consult \cite{Pet} or \cite{convex} for the necessary definitions.  

First of all, observe that if a quasigeodesic $\g:I\to X$ is horizontal at one point, it does not necessarily stay horizontal at every point. For example, consider the submetry from the square $R=[0,1]\times [0,1]$ onto $[0,1]$ given by projection onto one factor; the boundary of $R$ is a quasigeodesic, but at the corners changes from horizontal to vertical and vice versa. On the other hand, some properties do hold in this more general setting: 

\begin{prop}
The image of a horizontal quasigeodesic is a quasigeodesic in the base $B$. 
\end{prop}

\begin{proof}
Let $\g:I\to X$ be a horizontal quasigeodesic; denote by $\bg=\pi\circ\g$ its projection, which is clearly parametrized by arc length if $\g$ is. We will use the characterization of quasigeodesics given on page 171 of \cite{convex}: 
$\bg$ is a quasigeodesic iff there is an inequality 
$$
\angle(\bal'(0), \bg^+(0))\geq \tilde{\angle_k}(|\bg(0)\bq|,|\bg(t)\bq|,t)
$$
\noindent for all small $t>0$ and $\bq\in B$; here, $\bal$ is a shortest geodesic between $\bg(0)$ and $\bq$, $\bg^+(0)$ is the right tangent vector of $\bg$ at $0$, and $\tilde{\angle_k}$ is the comparison angle at $\bg(0)$ in a space of constant curvature $k$, where $k$ is a lower curvature bound for $X$ and $B$.  Since $\bal$ is minimal,  it may be horizontally lifted to a geodesic $\alpha$ starting at $\g(0)$ because of lemma \ref{lifts}.

We claim that the angle between the lifts equals the original angle; i.e, for horizontal $x$, $y$, $\angle(x,y)=\angle(D\pi x, D\pi y)$. To see this, recall that by Proposition \ref{hover}, there exists a minimal geodesic $\gamma:[0,1]\to H$ in $H$ joining $x$ and $y$, and the set 
$\{t\gamma(s)\mid t\ge0, s\in[0,1]\}$ in the cone $C(H)$ is an isometrically imbedded flat wedge. We may then define $x+y$ to be the midpoint of the segment joining $2x$ and $2y$. Even though this point depends on the choice of $\gamma$, one always has that $|x+y|^2=|x|^2+|y|^2+ 2\langle x, y\rangle$. Now, by a result of Lang and Schroeder \cite{LS97}, for any $v\in C(X)$, $\langle v, x+y\rangle \le \langle v,x\rangle + \langle v,y\rangle$. But by Proposition \ref{hover}, $w$ is horizontal if and only if $\langle v,w\rangle \le 0$ for all $v\in V$, so that $x+y$ is horizontal whenever $x$ and $y$ are. Finally, 
the angle between $x$ and $y$ has as cosine $(|x+y|^2-|x|^2-|y|^2)/(2|x||y|)$, and since all lengths are preserved under $D\pi$, it equals the angle between the projected vectors, as claimed.

Thus, in $X$, 
$$
\angle(\al'(0), \g^+(0))=\angle(\bal'(0), \bg^+(0)), \qquad |\g(t)q|\geq |\bg(t)\bq|,
$$ 
 where the last inequality is due to $\pi$ being distance nonincreasing. Hence 
$$
\tilde{\angle_k}(|\g(0) q|,|\g(t) q|,t) \geq \tilde{\angle_k}(|\bg(0)\bq|,|\bg(t)\bq|,t),
$$
and the result follows since 
$$
\angle(\al'(0), \g^+(0))\geq \tilde{\angle_k}(|\g(0) q|,|\bg(t)q|,t)
$$
because $\g$ is a quasigeodesic. 

\end{proof}


%
%
%

\section{Extremal subsets and submetries}\label{extremals}

Many of the differences between submersions and submetries seem to arise from the presence of ``singular'' sets in Alexandrov spaces. The appropriate version of singularity in this context is that of \emph{extremal set}, introduced by Perelman and Petrunin in \cite{PePe}. In this section, we examine how such sets are situated in relation to the fibers of the submetry. 
We refer the reader to \cite{convex} and the above reference for the definitions and lemmas used in this section: 


\subsection{Images of extremals are extremals}\label{image}
We generalize Proposition 4.1 in \cite{PePe} from isometric quotients to general submetries. 
\begin{prop}\label{P:extr}
Let $E\subset X$ an extremal set. Then $\bE:=\pi(E)$ is extremal in $B$. 
\end{prop}

\begin{proof}
Since $B$ itself is extremal, we may assume that $\bE$ is a proper subset of $B$. Let $\bq$ be a point in $B$, and suppose that $\bp$ is a point in $\bE$ at minimal distance from $\bq$.  Consider any $p\in E$ in the preimage of $\bp$, and choose some point $q$ over $\bq$ such that $|pq|=|\bp\bq|$. Then the distance function $\dis(q,\cdot)$ from $q$ has a minimum at $p$ when restricted to $E$, and for any sequence $p_i$ in $E$ converging to $ p$, 
$$
\limsup_{p_i\to p}\frac{|qp_i|-|qp|}{|pp_i|}\leq 0
$$
by the definition of extremality.
If $\bp_i$ is a sequence of points in $\bE$ approaching $\bp$, then after choosing $p_i$ over $\bp_i$ approaching $p$, 
using that $|\bq\bp_i|\leq |q p_i| $ and the above inequality, we get
$$
\limsup_{\bp_i\to \bp}\frac{|\bq\bp_i|-|\bq\bp|}{|\bp\bp_i|}
\leq \limsup_{p_i\to p}\frac{|qp_i|-|qp|}{|pp_i|}\cdot\frac{|pp_i|}{|\bp\bp_i|}
\leq 0.
$$  

\end{proof}

Observe that the converse is not true: if $\bar{E}\subset B$ is extremal, it does not, in general, follow that $\pi^{-1}(\bar{E})$ has that same property. For instance, take the action of $\Z_k$ on $\R^2$ where the generator of $\Z_k$ acts by a rotation with angle $2\pi/k$; if $k>1$, then the orbit space $B$ is a cone of angle $2\pi/k$, and thus its vertex $\bar{0}$ is extremal.  $\R^2$, however, is extremal-free.

\subsection{There are no slanted extremals} Recall that for a submetry $\pi:X\to B$, directions at points of $X$ split into horizontal and vertical components. Recall also that given any point $p$ in an extremal set $E\subset X$, there is a well defined subset $\sigE\subset\Sigma_p$ of directions tangent to $E$ at $p$. A very natural question arises: how is $\sigE$ placed with regard to $\pi$? 

It is easy to construct examples where $\sigE$ is entirely horizontal or vertical: 
\begin{itemize}
\item Let $K$ be a Euclidean cone over any positively curved Riemannian manifold of diameter less than $\pi/2$,  and define $X=K\times B$, for $B$ any Riemannian manifold. The only extremal set in $X$ is the product of the vertex of $K$ with $B$. It is clear that the projection onto $B$ is a submetry and $E$ is horizontal.
\item For the same $X$, choose now the projection onto $K$ to be the submetry: $E$ is now a fiber, hence vertical.
\item The above cases do not cover all the possibilities for extremals. After all, the total space $X$ is itself an extremal. But the situation can even mix horizontal and vertical directions without any intermediate behaviour. For instance, let $X$ denote the product $ C_0(\CP^2)\times  C_0(\CP^2)$ of the Euclidean cone over $\CP^2$ with itself, and $\pi$ the projection onto the first factor. Then the set $( C_0(\CP^2)\times\{0\})\cup (\{0\}\times  C_0(\CP^2))$ is extremal, but is neither horizontal nor vertical for the submetry $\pi$.
\end{itemize}

The following lemma is a particular case of our main structure result, but since it illustrates its proof, we state it separately:

\begin{lem}\label{suspension}
Let $X$be a spherical suspension over an Alexandrov space $Y$ with curvature $\geq 1$. If $E$ is a connected extremal set in $X$, then $E$ is either one of the cone points or a spherical suspension over some extremal subset of $Y$.
\end{lem}

\begin{proof}
Assume $E$ is not a cone point; we need to show that for any other $p\in E$, the geodesic $\g:[0,\pi]\to X$ connecting the poles through $p$ is contained in $E$. Recall that according to \cite{convex}, a set is extremal iff for any $q\in X$, a gradient curve of the distance function $\dis(q,\cdot)$ from $q$ starting at a  point of $E$ remains in $E$. But  the gradient curves of the distance functions from the poles are precisely the meridians. Thus we obtain that $\g[0,\pi]\subset E$ as desired. Denote by $A$ the set consisting of the midpoints $\g(\pi/2)$ of all meridians $\g$ that pass through points of $E$; then $A$ is contained in $Y$, which is totally geodesic. By looking at gradient curves from points in $Y$ we conclude that $A$ has to be extremal in $Y$ as claimed.  
\end{proof}

\begin{thm}[Structure of extremal directions]\label{extremal_tangent}
Let $\pi:X\to B$, $E$, $\sigE$ be as above, and set $HE = H\cap \Sigma_pE$, $VE = V\cap \Sigma_pE$.
\begin{enumerate}
\item  If $HE=\emptyset$, then $\sigE$ is vertical; an analogous statement holds for $VE$.
\item $HE$ and $VE$ are extremal in $H$ and $V$ respectively.
\item There is an inclusion $\sigE\subset P(HE,VE)$. 
\item Each connected component of $\sigE$ is either horizontal, vertical, or contained in $P(H_1,V_1)$ for some connected extremal sets $H_1$, $V_1$ in $H$ and $V$ respectively.
\end{enumerate}
\end{thm}


%

\begin{proof}
The starting point is that if $E$ is extremal in an Alexandrov space, then its space of directions $\sigE$ is extremal in $\sig_p$. 
Another useful fact is that  if $e=x+u$ is the horizontal-vertical decomposition of $e\in \sigE$, then the geodesic in $\sig$ connecting 
$x/|x|$ to $u/|u|$  must be contained in $\sigE$ whenever $x,u\neq 0$, since such a geodesic is a gradient curve for the distance function from $x/|x|$ or $u/|u|$ respectively. Combining these with the splitting of $\sig_p$ into horizontal and vertical directions implies most of the statements. 
 

  For instance,  choose some $e\in\sigE$ with $e=h+u$ its decomposition. If $HE=\emptyset$, then $h$ must vanish, since otherwise the observation in the first paragraph shows that $h/|h|\in \sigE$. Thus $e\in VE$. The case with $VE=\emptyset$ is proved in a  similar way. 

Assume now that $HE\neq \emptyset$, and let $h\in H$. Let $e=y+u\in\sigE$ be the point in $\sigE$ closest to $h$.  Then either $e=u$ or $e=y$ since otherwise by the above observation, $y$ would be a point in $\sigE$ closer to $h$ than $e$. If there were a point $z$  in $H$ at distance larger than $\pi/2$ from $h$, we would get a contradiction: a geodesic from $u$ to $h$ forms an angle 
$\dis_H(h,z)$ with a geodesic from $u$ to $z$, negating the extremality of $E$. Therefore the closest point to any $h\in H$ can be taken at $H$. Since this set is totally convex, $HE$ is extremal in $H$. Once again, the argument for $VE$ is entirely similar. 


Next assume that $w\in \sigE$. Write $w=h+v$; if $h$ or $v$ vanish, then clearly $w\in P(HE,VE)$. We may therefore assume that $w$ lies in the interior of some geodesic between $H$ and $V$.  As before, such a geodesic  remains in $\sigE$  and therefore $\sigE\subset P(HE,VE)$. Finally the last statement is obvious, since each connected component of an extremal set is itself extremal. 

\end{proof}

\begin{thm}
Let $E\subset X$ an extremal set, and $\pi:X\to B$ a submetry. Then $E$ is not slanted with respect to $\pi$.
\end{thm}
\begin{proof}
This is clear from the previous theorem since the argument shows that any point in $E$ must admit tangent directions that are either horizontal or vertical. 
\end{proof}

\begin{cor}
Let $E$ be a connected extremal set of dimension one in an Alexandrov space $X$. Then at each of its points $E$ is either horizontal or vertical with respect to any submetry from $X$.
\end{cor}
\begin{proof}
Otherwise $E$ would have some point $p$ with a slanted direction, and hence the above theorem would force $E$ to have dimension two or greater.
\end{proof}
In fact, as stated earlier there are one-dimensional extremal sets with both horizontal and vertical pieces. An easy argument shows that they are countable at most.

\subsection{The restriction of a submetry to an extremal set} 


Recall that extremal sets with their induced metrics are no longer  Alexandrov spaces in general. In spite of this,  they are nonetheless well behaved with respect to submetries. 
We make this more explicit in the following result, which is a particular case of the more general theorem 4.3 in \cite{Lyt_open}:

\begin{thm}
Let $E$ be an 
extremal subset of an Alexandrov space $X$, and $\pi:X\to B$ a submetry. 
Assume that for every $p\in E$ and for every direction $\bv\in\sig_{\pi(p)}\bE$ there is a $v\in\sigE_p$ with $D_p\pi(v)=\bv$. 
Then $\pi|_E:E\to \bE$ is a submetry when $E$ and $\bE$ are given their induced metrics.
\end{thm}

\begin{proof}
Let $p\in E$, $r>0$. If $\g:[0,r]\to E$ is a curve in $E$ with $\ell(\g)\leq r$, then $\ell(\pi\circ\g)\leq r$ and $\pi\circ\g$  is contained in $\bE$.  Hence $\pi(B^E(p,r))\subset B^{\bE}(\bp,r)$ where the balls are being taken in the relative metrics. 


For the opposite inclusion, let $\bq\in \bE$ be some point with $\dis_{\bE}(\bp,\bq)=r$; we want to show that there is some unit speed curve $\g:[0,r]\to X$ entirely contained in $E$, and such that 
$\g(0)=p$, $\g(r)\in\pi^{-1}(\bq)$. This will be constructed using piecewise quasigeodesics whose limit will result in the desired curve. 
To accomplish this, observe that the hypothesis on $D\pi$ combined with Theorem \ref{extremal_tangent} shows that any tangent direction to $\bE$ at some point $\bp'$ can be lifted horizontally and tangential to any point $p'$ in $E$ over $\bp'$. 

Given a positive integer $n$, consider the partition $0=t_0<\dots<t_i=\frac{i}{n}r<\dots< t_n=r$.  We construct $\g_n$ as follows: 
\begin{itemize}
\item Let $\g_n(0)=p$; lift a unit speed shortest $\bE$-curve between $\bp$ and $\bq$ to $p$, and denote it by $\al$. By Theorem \ref{extremal_tangent}, $\al'_+(0)$ is tangent to $E$; 
\item On the interval $[t_0,t_1]$, let $\g_n$ coincide with the (necessarily unit speed) 
quasigeodesic from $p$ tangent to $\al'_+(0)$; by extremality of $E$, this curve can be taken entirely contained in $E$
(see for instance section 5.1 in \cite{convex});
\item To extend $\g_n$ to $[t_1,t_2]$, take a $\bE$-geodesic in $B$ between $\pi\circ\g(t_1)$ and $\bq$, and  lift it horizontally to $\g_n(t_1)$. Again, its right velocity vector is tangent to $E$ by Theorem \ref{extremal_tangent}, and  we may follow a quasigeodesic from $\g_n(t_1)$ and contained in $E$  to construct $\g_n:[t_1,t_2]\to X$;
\item Iterate this procedure on each subinterval of the partition; after $n$ steps we obtain a curve $\g_n$ defined over all the interval $[0,r]$. 
\end{itemize}

There is a subsequence of the  $\g_n$(which we still denote by $\g_n$) that converges to some $\g:[0,r]\to X$. Clearly $\g$  is entirely contained in $E$, since extremal sets are closed. Thus, there are only two facts that remain to be established: firstly that the length of $\g$ does not exceed $r$, and secondly that $\g$ connects $p$ to the fiber over $\bq$. For the first one, observe that the length of each $\g_n$ does not exceed $r$, so that $\ell(\g)\leq r$. 

In order to prove that $\g(r)$ is in $\pi^{-1}(\bq)$, set $\bg=\pi\circ\g$, and observe that  we have  $\dis(\g(t), \pi^{-1}(\bq))\leq\dis_{\bE}(\bg(t),\bq)$. 
We will use Petrunin's  first variation formula for extremal sets (see \cite{Pet}), in order to show that the right hand side equals zero when $t=r$. So let $\bg_n=\pi\circ\g_n$; since at each $t_i$, $\bg_n'$ is tangent to a shortest geodesic from $\bg_n(t_i)$ to $\bq$, we have that for $t\in[t_i,t_{i+1}]$,
$$
|\bg_n(t)\bq|_E=|\bg_n(t_i)\bq|_E - (t-t_i)+o(t-t_i), 
$$ 
and in particular 
$$
|\bg_n(t_{i+1})\bq|_E=|\bg_n(t_i)\bq|_E - \frac{r}{n}+o(\frac{1}{n}).
$$
Writing out the above for $i=0,\dots,n-1$, we obtain
$$
|\bg_n(r)\bq|_E = |\bg_n(0)\bq|_E-r+n\cdot o(\frac{1}{n}) = n\cdot o(\frac{1}{n}).
$$
Taking the limit as $n\to\infty$ now implies that $\bg(r)=\bq$, as desired. 

\end{proof}

The hypothesis on $D\pi$ is necessary: otherwise let  $X=[0,1] \times[0,1]$, $B=[0,1]$, and $\pi$ be the projection onto the first factor. Then the boundary $E$ of $X$ is extremal, but at the point  $(0,1/2)$, say, balls of radius $r<1/2$ in $E$ are mapped to the point $0$.  

\section{A  splitting theorem for bounded curvature}

Recall that for length spaces $(X,\dis_X)$ and $(Y,\dis_Y)$, the \emph{product metric} $\dis$ on $X\times Y$ is defined by
\[\dis^2((x_1,y_1),(x_2,y_2))= \dis^2_X(x_1,x_2) + \dis^2_Y(y_1,y_2),\qquad x_i\in X, \quad y_i\in Y.
\]
One feature of the product metric is that the projections $\pi_X:X\times Y\to X$ and $\pi_Y:X\times Y\to Y$ are submetries for which the horizontal space of one is the vertical space of the other: formally, for each $(x,y)\in X\times Y$, $\pi_X:\pi_Y^{-1}(\pi_Y(x,y))\to X$ and $\pi_Y:\pi_X^{-1}(\pi_X(x,y))\to Y$ are isometries. For arbitrary length spaces, this does not characterize the product metric, however:
\begin{ex} Consider the metric $\tilde\dis$ on $X\times Y$, where
\[\tilde\dis((x_1,y_1),(x_2,y_2)) = \dis_X(x_1,x_2) + \dis_Y(y_1,y_2).
\]
It is easily checked that the projections are submetries. For example, given $x\in X$, $y\in Y$, if $B$ denotes the closed ball of radius $r$ in $X$ centered at $x$, and $\tilde B$ the closed ball of the same radius in $X\times Y$ centered at $(x,y)$, then $B\times\{y\}\subset \tilde B$, so that $B\subset\pi_X(\tilde B)$. Since $\pi_X$ does not increase distances, $\pi_X(\tilde B) = B$, and $\pi_X$ is a submetry. Clearly, the horizontal space at a point is the fiber of $\pi_Y$ through that point.
\end{ex}

For metrics with curvature bounded both above and below, however,  this property does characterize a product metric:
\begin{prop}
Let $\dis$ denote a metric on $X\times Y$ with curvature bounded above and below. Then $\dis$ is the product metric iff the projections onto each factor are submetries such that for each $(x,y)\in X\times Y$, $\pi_X:\pi_Y^{-1}(\pi_Y(x,y))\to X$ and $\pi_Y:\pi_X^{-1}(\pi_X(x,y))\to Y$ are isometries.
\end{prop}
\begin{proof}
The condition is clearly necessary, so assume that $\dis$ is a metric for which the projections are submetries. By a result of Nikolaev (see e.g.~\cite{BBI} Theorem 10.10.13), boundedness of curvature implies that $X\times Y$ is a $C^3$-manifold and $\dis$ is the distance function associated to a Riemannian metric  $\ip {}{}$ of class $C^{1,\al}$ for all $\al<1$. Identifying $T_{(x,y)}(X\times Y)\cong T_xX\times T_yY$ via $(\pi_{X*},\pi_{Y*})$, we have  $|(\mathbf{u},0)| = |\mathbf{u}|$ for $\mathbf{u}\in T_xX$ since $ (\mathbf{u},0)$ is $\pi_X$-horizontal, and similarly $|(0,\mathbf{v}))| = |\mathbf{v}|$ for $\mathbf{v}\in T_yY$. Furthermore,
\ip{(\bfu,0)}{(0,\bfv)}=0 by orthogonality of horizontal and vertical subspaces. Thus, 
\[|(\bfu,\bfv)|^2= \ip{(\bfu,0)+(0,\bfv)}{(\bfu,0)+(0,\bfv)}=|(\bfu,0)|^2+|(0,\bfv)|^2 +2\ip{(\bfu,0)}{(0,\bfv)} = |\bfu|^2+|\bfv|^2,
\]
which is the definition of the Riemannian product metric.
\end{proof}

Similar problems occur when attempting to generalize other  results for Riemannian submersions to the submetry setting: for example, it is known that the only complete warped products of nonnegative curvature are metric products. This fails for Alexandrov spaces whose curvature is not bounded above, as can be seen by considering a sharp cone $[0,\infty)\times _f S^1$, where $f(t) = \al t$, $0<\al<1$.

Another property of Riemannian submersions that does not carry over to submetries is uniqueness of horizontal lifts, as was pointed out in Section 2. For yet another example, consider the ray in the base space of the submetry  $[0,\infty)\times _f S^1\to[0,\infty)$ from above. It has infinitely many lifts at the vertex of the cone.

\end{document}